\documentclass[12pt,a4paper]{article}
\usepackage{amssymb}

\usepackage{graphicx,amssymb,amsfonts,epsfig,amsthm,a4,amsmath,url}
\usepackage[latin1]{inputenc}

\newtheorem{thm}{Theorem}[section]
\newtheorem{cor}[thm]{Corollary}
\newtheorem{lem}[thm]{Lemma}

\newtheorem{prop}[thm]{Proposition}
\theoremstyle{definition}
\newtheorem{defn}[thm]{Definition}

\theoremstyle{remark}

\theoremstyle{cremark}

\numberwithin{equation}{section}

\newcommand{\Z}{\mathbf{Z}}

\newcommand{\N}{\mathbf{N}}
\newcommand{\R}{\mathbf{R}}

\newcommand{\Diam}{\text{Diam}}

\newcommand{\bpr}{\noindent \textbf{Proof}: ~}

\newcommand{\epr}{~$\blacksquare$}

\newcommand{\eps}{\varepsilon}

\newcommand{\supp}{\textnormal{Supp}}

\newcommand{\A}{\mathcal{A}}

\title{Metric sparsification and operator norm localization}
\author{Xiaoman Chen, Romain Tessera, Xianjin Wang, Guoliang Yu}
\date{\today}

\begin{document}

\baselineskip=16pt

\maketitle

\begin{abstract}
We study an operator norm localization property and its applications to the coarse Novikov conjecture in operator K-theory.
A metric space $X$ is said to have operator norm localization property if there exists $0<c\leq 1$ such that for every $r>0$, there is   $R>0$
for which, if $\nu$ is a positive locally finite
Borel measure  on $X$,  $H$ is a separable infinite dimensional Hilbert space  and  $T$ is a bounded linear operator
acting on $L^2(X, \nu)\otimes H$ with propagation $r$, then there exists
an unit vector $\xi \in L^2(X, \nu)\otimes H$ satisfying the  $\Diam(\supp (\xi)) \leq R$ and $\| T\xi \| \geq c\|T\|.$
If $X$ has finite asymptotic dimension, then $X$ has operator norm localization property. In this paper, we introduce a sufficient geometric condition for the operator norm localization property. This is used to give many examples of finitely generated groups with infinite
asymptotic dimension and the operator norm localization property. We also show that any sequence of expanding graphs does not possess the operator norm localization property.

\end{abstract}


\section{Introduction}

Operator norm is a global invariant and is  often difficult to estimate. In this paper, we study a localization
property which allows us to estimate the operator norm locally relative to a metric space.
This property is motivated by the coarse Novikov conjecture
in operator K-theory.
We introduce a natural coarse geometric property on metric spaces, called
Metric Sparsification, to study the operator norm localization property. Roughly
speaking, this property says that  there exists a  constant $0<c\leq 1$ such that,
for every positive  finite Borel measure $\mu$
on $X$, there exists a subset $E$, which is a union of ``well
separated" subsets of ``controlled" diameters such that
$\mu(E)\geq c\mu(X)$. We show that the supremum over all $c$,
called the metric sparsification number of $X$, and denoted by $a(X)$ is a
coarse invariant. We show that the metric sparsification property implies the operator norm localization property, but is
more flexible than the latter. We prove for instance that any
solvable locally compact group equipped with a proper, locally
finite left-invariant metric has an approximation number $a(X)=1$.
This provides the first examples of finitely generated groups (as metric spaces with word metric) with infinite
asymptotic dimension satisfying operator norm localization property, as for instance
$asdim(\Z\wr\Z)=\infty$ ($\Z\wr\Z$ is the wreath product of $\Z$ with $\Z$). This also implies that connected Lie groups and
their discrete subgroups have a sparsification number equal to $1$. We obtain several permanence   properties for the operator norm localization property.  We  also
show that any sequence of expanding graphs does not possess operator norm localization property.
Finally in the last section of this paper, we apply the operator norm localization property to prove the coarse Novikov
conjecture for certain sequences of expanders.

The first, third and fourth authors of this paper are supported by CNSF. The second and fourth authors are supported
in part by NSF. The second and fourth authors wish to thank Fudan University for hospitality and support
during a visit in the summer of 2007.

\section{Operator norm localization}

In this section, we introduce an operator norm localization property for a metric space.
We show that this property is invariant under coarse geometric equivalence.

Recall that a Borel measure on a metric space is said to be locally finite if every bounded Borel subset
has finite measure.

\begin{defn}\label{defn} (Roe \cite{Roe1})
Let $X$ be a metric space with a positive locally finite Borel measure $\nu$, let $H$ be a separable
and infinite dimensional Hilbert space.
A bounded operator $T: L^2(X,\nu)\otimes H \rightarrow L^2(X, \nu)\otimes H,$ is said to have propagation
 at most $r$ if for all $\varphi,\psi\in L^2(X, \nu)\otimes H$ such that
$d\left(\supp(\varphi),\supp(\psi)\right)>r$,
$$\langle A\varphi,\psi\rangle=0.$$
\end{defn}

Note that if $X$ is discrete, then we can write
$$L^2(X,\nu)\otimes H =\oplus_{x\in X} (\delta_x \otimes H),$$ where $\delta_x$ is the Dirac function at $x$. Every bounded operator acting on $L^2(X,\nu)\otimes H$ has a corresponding matrix representation $$T=(T_{x,y})_{x,y\in X},$$ where
$T_{x,y}:\delta_y\otimes H  \rightarrow \delta_x\otimes H$ is a bounded operator. For $T$ to have propagation $r$, it
 is equivalent to saying that
the matrix coefficient $T_{x,y}$ of $T$ vanishes when $d(x,y)>r$.
The space of operators acting on $L^2 (X, \nu)\otimes H$ with propagation at most $r$ will be denoted by
$\mathcal{A}_r(X,\nu)$.

Let $\|T\|$ denote the operator norm of a bounded linear operator $T$.

\begin{defn}
Let $(X,\nu)$ be a metric space equipped with a positive locally finite Borel measure
$\nu$. Let $f$ be  a
(non-decreasing) function $\N\to \N$.  We say that $(X, \nu)$ has operator norm localization property relative to $f$  with constant $0<c\leq 1$ if, for all $k\in
\N$, and every $T\in \mathcal{A}_k(X,\nu)$, there exists nonzero
$\varphi\in L^2(X,\nu)\otimes H$ satisfying
\begin{itemize}
\item[(i)] $\Diam(\supp(\varphi))\leq f(k)$, \item[(ii)]
 $\|T\varphi\|\geq
c\|T\|\|\varphi\|$.
\end{itemize}
The supremum over all possible $c$ is called the operator
norm localization number of $(X,\nu)$ and is denoted by $opa(X,\nu).$
\end{defn}

\begin{defn}
A metric space $X$ is said to have operator norm localization property if  there exist a constant $0<c\leq 1$ and a
(non-decreasing) function $f:\N\to \N$ such that,
for every positive locally finite Borel measure $\nu$ on $X$, $(X, \nu)$ has
operator norm localization property relative to $f$ with constant $c$.
The supremum over all possible $c$ is called the operator
norm localization number of $X$ and is denoted by $opa(X).$
\end{defn}

We point out that a locally compact metric space $X$ has operator norm localization property if $(X, \nu_0)$
has operator norm localization property for some positive  locally finite Borel measure $\nu_0$ such that there exists $r_0>0$
for which every closed ball with radius $r_0$ has positive measure. This can be seen as follows.
We can decompose $X$ into countable  disjoint union of uniformly bounded Borel subsets $\{X_i\}_{i\in I}$
such that every bounded subset of $X$ is contained in a union of finitely many members of  $\{X_i\}_{i\in I}$ and
 $\nu_0 (X_i)>0$. We decompose
$$L^2(X,\nu_0)\otimes H=\oplus_{i\in I} (L^2(X_i, \nu_0)\otimes H).$$
For every other positive locally finite Borel  measure $\nu$, we have a similar decomposition:
$$L^2(X,\nu)\otimes H=\oplus_{i\in I} (L^2(X_i, \nu)\otimes H).$$
Let $W: L^2(X,\nu)\otimes H \rightarrow L^2(X,\nu_0)\otimes H,$ be an isometry such that
$$W (L^2(X_i, \nu)\otimes H)\subseteq  L^2(X_i, \nu_0)\otimes H$$
for every $i\in I$. If $T$ is a bounded operator acting on $L^2(X,\nu)\otimes H$ with propagation $r>0$, then
$WTW^*$ is a bounded operator acting on $L^2(X,\nu_0)\otimes H$ with propagation $r+2D$ and $\|WTW^*\|=\|T\|$, where $D=sup\{Diam(X_i):i\in I\}$.
It follows that if $(X, \nu_0)$ has operator norm localization property relative to $f$ with constant $0<c\leq 1$, then
$(X, \nu)$ has operator norm localization property relative to $f+D$ with constant $c$.

Let  $F$ be  a Borel  map from a metric space $X$ to another metric space $Y$.
Recall that $F$ is said to be coarse if (1) for every $r>0$, there exists $R>0$
such that $d(F(x), F(y))<R$ for every pair of points $x$ and $y$ in $X$ satisfying
$d(x,y)<r$; (2) the inverse image $F^{-1}(B)$  for every bounded subset $B$ of $Y$ is bounded.
We say that $X$ is coarse equivalent to $Y$ if there exist coarse maps $F: X\rightarrow Y$
and $G:Y\rightarrow X$, such that there exist a constant $C$ satisfying $d(G(F(x)), x)<C$ for
all $x\in X$, and $d(F(G(y)), y))< C$ for all $y\in Y.$

\begin{prop}
The operator norm localization property is invariant under coarse equivalence. More precisely, if $X$ is coarse equivalent
to $Y$, then $opa(X)=opa(Y)$.
\end{prop}
\bpr
Let $X$ and $Y$ be two coarse equivalent metric spaces. Let $F: X\rightarrow Y$ and $G: Y \rightarrow X$
be two coarse maps as in the definition of coarse equivalence.
 There exist two increasing
functions $\rho_1,\rho_2:\R_+\to\R_+$ such that $\lim_{t\to
\infty}\rho_1(t)=\infty$ and $$\rho_1(d(x,x'))\leq
d(F(x),F(y))\leq \rho_2(d(x,x')).$$  We shall prove that if $Y$ has operator norm localization property with constant $0<c\leq 1$,
then so does $X$.

 Let
$\nu$ be a  positive locally finite measure on $X$ and let
$\nu'=F(\nu)$. It is not difficult to see that there exists an isometry $$W: L^2(X,\nu)\otimes H\rightarrow L^2(Y,\nu')\otimes H$$ satisfying $$\supp(W\varphi)\subseteq \{y\in Y: d(y, F(\supp (\varphi)))\leq 1\}.$$ For every $T \in \mathcal{A}_k(X,\nu)$,
we have that $\|T\|=\|WTW^*\|$ and $WTW^*\in \mathcal{A}_{k+1}(Y,\nu')$. These properties of $W$ imply that
$X$ has operator norm localization property.
\epr

Recall that a metric space $X$ is said to coarsely embed into $Y$ if $X$ is coarse equivalent to a subset of $Y$ (with the metric
induced from $Y$). The proof of Proposition 2.4 shows the following:

\begin{prop} If a metric space $X$ coarsely embeds into another metric space $Y$, then $opa(X)\geq opa(Y)$.
\end{prop}

It is an open question to find a geometric condition equivalent to the operator norm localization property.

\section{Metric sparsification property  and sparsification number}

In this section, we introduce the metric sparsification property and the sparsification number.
We prove that the sparsification number is a coarse geometric  invariant. In particular,
we show that any locally compact solvable group has sparsification number 1. As a consequence,
every connected Lie group has sparsification number 1.

\begin{defn}\label{Defn}
Let $X$ be a metric space. We say that $X$ has Metric Sparsification Property
 with constant $0<c\leq 1$ (for short we say that $X$ has
MS(c)), if there exists a (non-decreasing) function $f:\N\to \N$
such that for all $m\in \N$, and every finite positive Borel
measure $\mu$ on $X$, there is a Borel subset $\Omega=\sqcup_{i\in
I}\Omega_i$ such that
\begin{itemize}
\item[(i)] $d(\Omega_i,\Omega_j)\geq m$ for all $i\neq j\in I$,
\item[(ii)] $\Diam(\Omega_i)\leq f(m)$ for all $i\in I$,
\item[(iii)]  $\mu(\Omega)\geq
c\mu(X)$.
\end{itemize}
The supremum over all possible $c$ is called the sparsification
number of $X$ and is denoted by $a(X).$
\end{defn}

If the asymptotic dimension of a metric space $X$ is $n$, then $a(X)$ is greater
then or equal to $\frac{1}{n+1}$.

When we need to be more explicit, we will say that $X$ has MS(c)
with function $f$. If $m$ and $\mu$ are given, and if we want to
say that a subset $\Omega$ satisfies the conditions of
Definition~\ref{Defn}, we will simply write
$\Omega=\Omega(\mu,f,m,c)$.
\begin{defn}
We say that a family of metric spaces has uniform MS(c) if there
is an $f$ that works for all the elements of the family.
\end{defn}

\begin{prop}\label{coarseProp} Let $X$ and $Y$ be two metric spaces.
If $F:X\to Y$ is a coarse embedding, then $a(X)\geq a(Y)$.
\end{prop}
\bpr As $F$ is a coarse embedding, there exist two increasing
functions $\rho_1,\rho_2:\R_+\to\R_+$ such that $\lim_{t\to
\infty}\rho_1(t)=\infty$ and $$\rho_1(d(x,x'))\leq
d(F(x),F(y))\leq \rho_2(d(x,x')).$$ Assume that $Y$ has MS(c). Let
$\mu$ be a finite measure on $X$ and let $m\in \N$. Let
$\mu'=F(\mu)$, and let $\Omega'=(\mu',f,m,c)$ (for some $f$). Let
$\Omega=F^{-1}(\Omega')$. Then one immediately checks that
$\Omega=\Omega(\mu,\rho_1^{-1}\circ f\circ \rho_2 ,
\rho_2^{-1}(m),c)$. \epr

\begin{cor}\label{coarseEqCor}
The sparsification number is invariant under coarse equivalence.
\end{cor}

\begin{cor}
Let $X'$ be a metric space and let $X\subset X'$ be a metric
subspace of $X'$, i.e. a Borel subset of $X'$ equipped with
the induced distance. Then $a(X)\geq a(X')$.
\end{cor}

Let us prove an easy but useful lemma.

\begin{lem}\label{lem1}
Let $X$ be a metric space and assume that for every $m\in \N$,
there is an $m$-disjoint family of metric subspaces $(X_j)_{j\in
I}$ with uniform MS(c). Then there is a function $f$ such that for
every $m\in \N$, and every finite measure $\mu$ on $X$ supported
on $\bigcup_j X_j$, there exists $\Omega=\Omega(\mu,f,m,c)$
(included in $\bigcup_j X_j$).
\end{lem}
\bpr Let $\mu_j$ be the restriction of $\mu$ to $X_j$. As $(X_j)$
has uniform MS(c), there is a function $f_m$ such that for every
$j\in J$, there exists $\Omega_j=\Omega_j(\mu_j,f_m,m,c)$. Now
take $\Omega=\bigcup_j \Omega_j$. Clearly,
$\Omega=\Omega(\mu,f,m,c)$ for $f(m)=f_m(m).$ \epr

\begin{prop}\label{opensubProp}
Let $G$ be a locally compact group such equipped with some proper,
locally bounded left-invariant metric $d$. Let $G_n$ be a
non-decreasing, exhaustive sequence of open subgroups of $G$. If
the $G_n$ have MS(c) for the same constant $c>0$, then so does
$G$.
\end{prop}
\bpr Say that for all $n\in \N$, $G_n$ has MS(c) with function
$f_n$. Fix $m\in \N$ and take a finite measure $\mu$ on $G$. As
$G$ is locally compact and $d$ is proper, there exits $n=n(m)$
such that $B(1,m)\subset G_n$. Hence  the set of left cosets of $G$ modulo $G_n$ is an
$m$-separated partition of $G$. Let $\mu_n$ be the restriction of
$\mu$ to $G_n$. Let $\Omega_n=\Omega_n(\mu_n,f_n,m,c)$. Then
$\Omega=\Omega_{n(m)}=\Omega(\mu,f,m,c)$ where
$f(m)=f_{n(m)}(m)$.\epr

\begin{prop}\label{extensionProp}
Let $G$ be a locally compact compactly generated group and let $N$
be a closed normal subgroup of $G$. Assume that $N$ has MS(c) for
the induced metric, and that $G/N$ has MS($c'$). Then $G$ has
MS($cc'$).
\end{prop}
\bpr To fix the ideas, we equip $G$ with the word metric
associated to a compact generating subset $S$, and $G/N$ with the
word metric associated to $\pi(S)$, where $\pi$ is the projection
$\pi:G\to G/N$.

Fix $m\in \N$ and take some finite measure $\mu$ on $G$. Let
$\overline{\mu}=\pi(\mu)$, and let
$\overline{\Omega}=\overline{\Omega}(\overline{\mu},\overline{f},m,c')$.
Hence $\overline{\Omega}=\sqcup_{i\in I}\overline{\Omega}_i$,
where $Diam(\overline{\Omega}_i)\leq \overline{f}(m)$ and
$d(\overline{\Omega}_i,\overline{\Omega}_j)\geq m$ for all $i\neq
j$. Write $X_i=X_i(m)=\pi^{-1}\left(\overline{\Omega}_i\right)$.
As $\pi$ is 1-Lipschitz, we have that for all $i\neq j$,
$$d\left(\pi^{-1}(X_i),\pi^{-1}(X_j)\right)\geq m.$$

For every $i\in I$, let $x_i\in \overline{\Omega}_i$ and let
$g_i\in G$ such that $x_i=\pi(g_i)$. As $N$ is normal and
$X_i\subset B_Q(x_i,m)=x_iB_Q(1,m)$, we have
\begin{eqnarray*}
g_iN \subset \pi^{-1}(X_i) &\subset &
g_iB_G\left(1,\overline{f}(m)\right)N \\
& = & g_iN
B_G\left(1,\overline{f}(m)\right)\\
& = & \{g\in G, d(g,x_iN)\leq \overline{f}(m)\}.
\end{eqnarray*}
Hence, the obvious injection $N\to X_i$ where $1$ is sent to $x_i$
is a coarse equivalence with $\rho_1$ and $\rho_2$  depending only on $m$, hence
uniform in $i$. Hence the family $X_i$
has uniform MS(c) (see the proof
of Proposition~\ref{coarseProp}). Let $\mu_m$ be the restriction to $\bigcup_i
X_i$. Note that $\mu_m(G)=\overline{\mu}(\overline{\Omega})\geq
c'\overline{\mu}(G/N)=c'\mu(G).$ By Lemma~\ref{lem1}, there exists
$\Omega=\Omega(\mu_m, f, m,c)$ (for some $f$). Hence, together
with the previous remark, it yields $\Omega=\Omega(\mu, f,
m,cc').$ \epr

\begin{thm}
The sparsification number of a solvable locally compact group
equipped with a proper locally finite left invariant metric equals
$1$.
\end{thm}

As a consequence, the sparsification number of the  wreath product $\Z\wr\Z$ is $1$.
This implies that $\Z\wr\Z$ has operator norm localization property despite the fact
it has infinite asymptotic dimension.

\bpr By proposition~\ref{opensubProp}, we can assume that $G$ is
compactly generated. By Proposition~\ref{extensionProp}, we can
assume that $G$ is abelian, and then again
Proposition~\ref{opensubProp} reduces the problem to $G=\Z$ or
$\R$, so by Proposition~\ref{coarseProp}, to $\Z$.

Fix $m\in \N$ and let $\mu$ be a probability on $\Z$. Let $k\in
\N$, and let $\pi$ be the projection $\Z\to \Z/m(k+1)\Z$. For
every $j=0,\ldots, k$, let $C_j=\pi^{-1}\left([jm,(j+1)m)\right).$
The $C_j$, $j=0,\ldots, k$ form a partition of $\Z$. Hence, there
exists $j_0$ such that
$$\mu(C_{j_0})\leq 1/k.$$ Thus the complement $\Omega$ of $C_{j_0}$ in
$\Z$ satisfies $\Omega=\Omega(\mu,f,m,1-1/k)$ with $f(m)=km$. \epr

\begin{cor}
Every Borel  subset of a connected Lie group has an
sparsification number equal to $1$.
\end{cor}
\bpr As every connected Lie group has a co-compact solvable closed
group, it is coarse equivalent to a solvable group. \epr

It is an open question whether every finitely generated linear group equipped with a word metric
has metric sparsification property or operator norm localization property. It is also an interesting question to
compute the sparsification number (or operator norm localization number) of a simply connected and non-positively curved Riemannian manifold.

\section{Link to operator norm localization property}

In this section, we show that the metric sparsification property implies the
operator norm localization property. More precisely, we have the following result.

\begin{prop}
Let $X$ be a metric space. Then
$$\sqrt{a(X)}\leq opa(X),$$ where $a(X)$ and $opa(X)$ are respectively the sparsification number and the operator norm localization number of $X$.
\end{prop}
\bpr Let $\nu$ be a
 positive locally finite Borel  measure  on
$X$.
 For any measurable subset $U$ of $X$, let $P_U$ be the
orthogonal projector on the space of functions of $L^2(X)\otimes H$
supported on $U.$ Clearly, $A\in \mathcal{A}_k$ means that for any
subsets $U,V\subset X$ such that $d(U,V)>k$, $P_U AP_V=0.$ We
deduce that  if $\psi\in L^2(X)\otimes H$ is supported in $U$, then $A\psi$
is supported in $[U]_k:=\{x\in X, d(x,U)\leq k\}$. As a result, we have
\begin{lem}\label{lem2}
If $\psi\in L^2(X,\nu)\otimes H$ is a sum of non-zero $\psi_i\in
L^2(X,\nu)\otimes H$ whose supports are piecewise at distance larger than
$m> 2k$, then $$d\left(\supp(A\psi_i),\supp(A\psi_j)\right)\geq
m-2k>0$$ for all $i\neq j$.

Consequently,
$$\frac{\|A\psi\|}{\|\psi\|}\leq \sup_{i\in
I}\frac{\|A\psi_i\|}{\|\psi_i\|}.\;$$
\end{lem}
Now let $k\in \N$ and $A\in \mathcal{A}_k(X,\nu)$. Let $\varphi\in
L^2(X,\nu)\otimes H$. Consider the finite measure $d\mu=\|A\varphi\|_H^2d\nu$
and some $m>2k$. Let $\Omega=\bigcup_{i\in
I}\Omega_i=\Omega(\mu,f,3m,c)$ for some $c<a(X)$, where $f$ is as in Definition 3.1. Let $P_{\Omega}$
be the orthogonal projector on $L^2(\Omega)\otimes H$. Therefore,
$P_{\Omega}\varphi$ is a sum of $\varphi_i=P_{\Omega_i}\varphi\in
L^2(X,\nu)\otimes H$ (which we can assume to be non-zero) whose supports
are piecewise at distance larger than $3m$ and have diameter at
most $f(3m)$. Let $[\Omega]_m=\{x\in X, d(x,\Omega)\leq m\}$. Note
that $\Omega'=[\Omega]_m=\bigcup_{i\in I}[\Omega_i]_m$, and that
$\Omega'=\Omega'(\mu,f',m,c)$ with $f'(m)=f(3m)$.
\begin{lem}\label{lem3}
For all $\psi\in L^2(X,\nu)\otimes H$, $\|AP_{[U]_m}\psi\|\geq
\|P_{U}A\psi\|.$
\end{lem}
\bpr As $A\in \A_k(X,\nu)$ and $m>2k$, we have $P_{U}A
P_{X\smallsetminus[U]_m}=0$. Hence $P_UAP_{[U]_m}=P_UA$. So the
lemma follows.\epr

Using the first part of Lemma~\ref{lem2} and Lemma~\ref{lem3}, we
obtain
\begin{eqnarray*}
\|AP_{\Omega'}\varphi\|^2 &=& \sum_i\|AP_{[\Omega_i]_m}\varphi\|^2\\
& \geq & \sum_i\|P_{\Omega_i}A\varphi\|^2\\
& = &\|P_{\Omega}A\varphi\|^2.
\end{eqnarray*}
Hence,
$$\frac{\|AP_{\Omega'}\varphi\|^2}{\|P_{\Omega'}\varphi\|^2}
\geq  \frac{\|AP_{\Omega'}\varphi\|^2}{\|\varphi\|^2}  \geq
\frac{\|P_{\Omega}A\varphi\|^2}{\|\varphi\|^2}=
\frac{\mu(\Omega)}{\|\varphi\|^2}\geq c
\frac{\mu(X)}{\|\varphi\|^2}=c
\frac{\|A\varphi\|^2}{\|\varphi\|^2}.$$ Applying this
inequality to some $\varphi\in L^2(X,\nu)\otimes H$ such that
$\|A\varphi\|/\|\varphi\|\geq (1-\eps)\|A\|$, and applying
the second part of Lemma~\ref{lem2} to $\psi=P_{\Omega'}\varphi$,
we get that $X$ has operator norm localization property with constant $\sqrt{c}(1-\eps )$, for arbitrary small
$\eps>0$.\epr

It is an open question whether the metric sparsification property is equivalent to
the operator norm localization property.

\section{Permanence properties for  operator norm localization}

In this section, we prove several permanence properties for the operator norm localization properties.

Let $\Gamma$ be a group acting on a metric space $X$. For every $k\geq 0$,  the $k$-stabilizer $W_k(x_0)$
of a point $x_0\in X$ is defined to be the set of all $g\in \Gamma$ with $gx_0\in B(x_0, k)$, where $B(x_0,k)$ is the closed ball
with center $x_0$ and radius $k$. The concept of $k$-stabilizer is introduced by Bell and Dranishnikov in their work
on permanence properties of asymptotic dimension \cite{BD}.

\begin{prop} Let $\Gamma$ be a finitely generated group acting isometrically on a metric space $X$. If $X$ has metric sparsification property with a constant $0<c\leq 1$ and there exist $0<c'\leq 1$ and $x_0\in X$ such that $W_k(x_0)$ has operator norm localization property with constant $c'$ for each $k>0$, then
$\Gamma$ has operator norm localization property with constant $\sqrt{c}c'$ as a metric space with a word metric.
\end{prop}

\bpr
We define a map $\pi: \Gamma \rightarrow X$ by: $\pi (g)= g x_0$ for all $g\in \Gamma$.
Let $S$ be the finite generating set in the definition of the word metric for $\Gamma$ and let $\lambda = max\{ d(\gamma x_0, x_0), \gamma\in S\}$. It is easy to see that $\pi$ is $\lambda$-Lipschitz, i.e. $d(\pi (x), \pi(y))\leq \lambda d(x,y)$ for all $x$ and $y$ in $\Gamma$.

Let $\nu$ be a positive locally finite measure on $\Gamma$ and $H$ be a separable infinite dimensional Hilbert space. Let
$T: \ell^2(\Gamma, \nu)\otimes H \rightarrow \ell^2(\Gamma, \nu)\otimes H,$ be a bounded linear operator with propagation $r$ for
some $r>0$. For each vector $v\in \ell^2(\Gamma, \nu)\otimes H,$ we define a finite measure $\mu$ on $\Gamma x_0$ by:
$$\mu(\{x\})=\|P_{\pi^{-1}(x)} Tv\|_{ \ell^2(\Gamma, \nu)\otimes H}^2$$
for every $x\in \Gamma x_0$, where $P_{\pi^{-1}(x)}$ is the projection from  $\ell^2(\Gamma, \nu)\otimes H$ to its subspace
$\ell^2(\pi^{-1}(x), \nu)\otimes H.$

By the definition of metric sparsification property,  there exists
a subset $\Omega=\sqcup_{i\in
I}\Omega_i$ of $\Gamma x_0$ such that
\begin{itemize}
\item[(i)] $d(\Omega_i,\Omega_j)\geq (\lambda +10)(r+10)$ for all $i\neq j\in I$,
\item[(ii)] $\Diam(\Omega_i)\leq D$ for some $D>0$ and all $i\in I$,  where $D$ is independent of $\nu$ and $v$,
\item[(iii)]  $\mu(\Omega)\geq c\mu( \Gamma x_0 )$.
\end{itemize}

Notice that there exists $k>0$ such that $\pi^{-1}( \Omega_i)$ is coarse equivalent to a subset of $W_k(x_0)$
for all $i\in I$ with a uniform $\rho_1$ and $\rho_2$, where $\rho_1$ and $\rho_2$ are control functions as
in the proof of Proposition 2.4. By Proposition 2.4,  $\pi^{-1}( \Omega_i)$ has uniform operator
norm localization property with constant $c'$ for all $i\in I$ in the sense that each $\Omega_i$ has operator norm localization property  with constant $c'$ and  the function $f$ in the Definition 2.3 is
independent of $i\in I$.  This, together with the above properties of $\Omega_i$ and the fact
that $T$ has propagation $r$, implies that $$\|P_{\pi^{-1}(\Omega)}Tv\|^2\geq c\|Tv\|^2$$ and $P_{\pi^{-1}(\Omega)}T$ decomposes $$P_{\pi^{-1}(\Omega)}T=\oplus_{i\in I} T_i,$$ where  each $T_i$ is
an operator   acting on $\ell^2(\{g\in \Gamma: d(g, \pi^{-1}(\Omega_i))\leq r\})\otimes H$ with propagation $r$.
Note that $\{g\in \Gamma: d(g, \pi^{-1}(\Omega_i))\leq r\}$ is uniformly coarse equivalent to $\pi^{-1}( \Omega_i)$ and hence has
uniform operator norm localization property with constant $c'$ for all $i\in I$.
It follows that $\Gamma$ has operator norm localization property with constant $\sqrt{c}c'$.
\epr

Next we shall prove the following countable union result for operator norm localization property.
We say that a family of metric spaces $\{X_i\}_{i\in I}$ has uniform operator norm localization property with constant $0<c\leq 1$
if $X_i$ has operator norm localization property with constant $c$ for each $i$ and  the function $f$ in the Definition 2.3 is
independent of $i\in I$.

\begin{prop} Let $X$ be a metric space and $X=\cup_{i\in I} X_i$, where each $X_i$ is a Borel subset of $X$.
If $\{X_i\}_{i\in I}$ has uniform operator norm localization property with constant $0<c\leq 1$
 and, for each $r>0$, there exists  a Borel subset $Y_r\subseteq X$ having operator norm localization property with constant $c$ such that $\{X_i-Y_r\}_{i\in I}$ is $r$-disjoint, then $X$ has
operator norm localization property with constant $c-\epsilon$ for every $\epsilon>0$.
\end{prop}

\bpr
Let $\nu$ be a positive locally finite Borel measure on $X$.
Let $T$ be a bounded  linear operator acting on $L^2(X,\nu)\otimes H$ with propagation $r>0$.
 For every $1>\delta>0$, there exists a
unit vector $\xi\in L^2(X, \nu)\otimes H$
satisfying $\| T\xi\|\geq (1-\delta)\|T\|.$

Let  $$Z_k=\{x\in X: 10(k-1)\leq   d(x, Y_{10r})< 10(k+1)r\}$$ for each $k\in\Bbb{N}$.
Let $\xi_k \in L^2(X, \nu)\otimes H $ be defined by:   $\xi_k (x)=\xi(x)$  for all $x\in Z_k$ and
$\xi_k (x)=0$ for all $x\in X-Z_k$. We have $\|\xi\|^2=\sum_{k} \|\xi_k\|^2$. Hence for each large  $N\in \Bbb{N}$, there exists
$k_0\in \Bbb{N}$ satisfying $\|\xi_{k_0}\| < \frac{1}{N}.$

Let $U_1= \cup_{k<k_0} Z_k$ and $U_2=\cup_{k>k_0}Z_k$. Notice that $U_1$ and $U_2$ are $10r$-disjoint if both $U_1$ and $U_2$
are non-empty. By our assumptions and Proposition 2.4, the $r$-neighborhood of $U_1$ has operator norm localization property with constant $c$
 and the $r$-neighborhood of $U_2$ is the  union of pairwise $5r$-disjoint subsets having  uniform operator norm localization property with constant $c$.
Let $P_i$ be the projection from $L^2(X, \nu)\otimes H$ onto $ L^2(U_i, \nu)\otimes H$ for $i=1, 2$.
By the choice of $k_0$ and the fact $T$ has propagation $r$, we have
$$\max\{\|TP_1\|, \|TP_2\| \}\leq \|T\| \leq (\frac{1}{1-\delta}+\frac{1}{N}) \max\{\|TP_1\|, \|TP_2\| \}.$$
The above inequality, together with  our assumptions and the fact that $TP_i$ has propagation $r$ and is supported on the $r$-neighborhood of $U_i$,   implies our result.
\epr

\begin{cor} Let $A$ and $B$ be two finitely generated groups with a common subgroup $C$. The amalgamated product
$A\ast_C B$ has operator norm localization property if and only $A$ and $B$ have operator norm localization properties.
\end{cor}

\bpr
It is enough to prove the ``if'' part. We follow the strategy in Bell and Dranishnikov \cite{BD}. Bell and Dranishnikov
constructed a tree on which  $A\ast_C B$ acts isometrically \cite{BD}. Recall that a tree has asymptotic dimension $1$ and hence has operator norm localization property.
Proposition 5.2, together with the argument in the proofs of  Theorem 5 and Proposition 4 in \cite{BD},
shows that the $k$-stabilizer of this action  has operator norm localization property for each $k>0$. Our corollary now follows from Proposition 5.1.
\epr

By using Propositions 5.1, 5.2 of this paper and constructions in section 5 of \cite{BD}, we can prove the following permanence result  for operator norm localization property in the case of  HNN extensions.

\begin{cor} Let $G$ be a finitely generated group with a word metric.
Let $\phi: A\rightarrow G,$ be a monomorphism of a subgroup $A$ of $G$, let $G'$ be the HNN extension of $G$.
If $G$ has  operator norm localization property, then $G'$ has operator norm localization property.
\end{cor}

We should point out that similar permanence results for finite asymptotic dimension was obtained by Bell and Dranishnikov in \cite{BD}.

\section{Expanding graphs and operator norm localization property}

In this section, we show that that any expanding sequence of  graphs doesn't have operator norm localization property.
In particular, this implies that any expanding sequence of graphs doesn't have the metric sparsification  property
defined in this paper.

For convenience of readers, we briefly recall the concept of expanding graphs \cite{Lub}.

\begin{defn}
Let $X=X(V, E)$ be a finite graph with $V$ as its vertex set and $E$ as its edge set.
Define the Cheeger constant of $X$ by:

$$ h(X) =\inf_{A, B\subseteq V} \frac{|E(A, B)|}{\min(|A|, |B|)},$$
where the infimum is taken over all disjoint partition $V= A\cup B$ and
$E(A, B)$ is the set of all edges connecting vertices in $A$ to vertices in $B$.

\end{defn}

\begin{defn}
An infinite  sequence of graphs $\{X_n(V_n, E_n)\}_{n=1}^{\infty}$ of bounded degree is said to be a
sequence of expanding graphs if there exists $h>0$ such that  $h(X_n)\geq h$ for all $n$ and the number of elements in $V_n$ goes to $\infty$ as $n\rightarrow \infty$.
\end{defn}

In the sense of probability, most sequences of graphs are expanding \cite{Lub}.

\begin{defn} The Laplacian $\triangle$ of the graph $X=X(V, E)$ is the operator on $l^2(V)$
defined by:

$$\triangle f(x)=\sum_{y\in V} \delta_{xy} (f(x)-f(y))$$
for every $f\in \ell^2(V)$,
where $\delta_{x,y}$ is the number of edges between $x$ and $y$,
\end{defn}

It is not difficult to show that $\triangle$ is self-adjoint and positive.
Let $\lambda_1(X)$ be the smallest positive eigenvalue.

The following result is well-known \cite{Lub}.

\begin{prop}
$\{ X_n\}_{n}$ is an expanding sequence of
 graphs if and only if there exists $\lambda>0$
such that $\lambda_1 (X_n)\geq \lambda$ for all $n$.
\end{prop}

Let $\{X_n\}_{n=1}^{\infty}$ be an infinite sequence of  graphs.
We endow a metric on the disjoint union $\cup_{n}X_n$ such that the restriction
of the metric on each connected component of $X_n$ is the natural path metric and $d(X_i, X_j)>i+j$ if $i\neq j$.
Let $V=\cup_n V_n \subseteq \cup_n X_n$ be given its subspace metric.

\begin{thm}
If $\{X_n\}_n$ is an infinite expanding sequence of  graphs, then
the metric space $V $ defined as above doesn't have operator norm localization property.
\end{thm}

\bpr
Let $\triangle_n$ be the Laplacian of the graph $X_n$.
Let $p_n$ be the projection from $\ell^2(V_n)$ to the one dimensional subspace
of constant functions. By abuse of notation, we denote the operator  $p_n\otimes I$ acting on $\ell^2 (V_n)\otimes H$
by $p_n$ and the operator $\triangle_n\otimes I$ acting on $\ell^2 (V_n)\otimes H$ by
$\triangle_n.$
 Let $p=\oplus_{n} p_n$ and $\triangle=\oplus_n \triangle_n$.
We have $$p=\lim_{t\to +\infty} exp(-t\triangle), $$
where the limit is taken in operator norm (as operators acting on the Hilbert space
$\ell^2(V)\otimes H$ ).
It follows that, for any $\epsilon>0,$ there exist an operator $T$  and $r>0$ in $B(\ell^2(V)\otimes H)$  such that
$||T-p||<\epsilon$ and $T$ has propagation $r$.
The fact that $T$ has finite propagation implies that there exists some large  $N$ such that
$\ell^2 (V_n)\otimes H$ is invariant under $T$ and $T^*$ if $n>N$. We denote the restriction of $T$ to
$\ell^2(V_n)\otimes H$ by $T_n$ if $n>N$. Consider
$$S_k=\oplus_{n\leq k} 0 \oplus_{n>k} T_n\in B(\ell^2(V)\otimes H)$$ if $k\geq N$
and
$$Q_k=\oplus_{n\leq k} 0 \oplus_{n>k} p_n\in B(\ell^2(V)\otimes H)$$ if $k\geq N$.
 We observe that $S_k$ has propagation $r$ and
$||S_k-Q_k||<\epsilon$ if $k\geq N$.

Now we assume by contradiction that $V$ has operator norm localization property. By assumption, there exist $C>0$ (independent of
$\epsilon$ and $r$) and $R$ (dependent on $r$) and a unit vector
$v_k\in \ell^2(V)\otimes H$
such that $||S_k||\leq C ||S_k v_k|| $ and $\Diam (\supp (v_k))< R$ for all $k>N$.
We have $$||S_kv_k||\leq ||Q_kv_k|| +\epsilon$$
for all $k>N$. By the definition of $Q_k$ and the support condition of $v_k$,
we know $||Q_kv_k||\rightarrow 0$ as $k\rightarrow \infty$.
Consequently we have $$||S_k||\leq (1+C)\epsilon$$ if $k$ is large enough.
However, by the definition of $S_k$ and the fact that $Q_k$ has norm $1$, we have $$||S_k||> 1-\epsilon$$ if $k>N$.
This is a contradiction if we choose $\epsilon$ small enough and $k$ large enough.
\epr

\section{Applications to K-theory}

In this section, we discuss applications of the operator norm localization property to the coarse Novikov conjecture.

Let $\Gamma$ be a finitely generated residually finite group. We can assume that there is a sequence of normal
subgroups of finite index
$$\Gamma_1\supseteq \Gamma_2\supseteq \cdots \supseteq \Gamma_i\supseteq \cdots$$
such that
$$\bigcap_{i=1}^\infty \Gamma_i =\{e\}.$$
Endow $\Gamma/\Gamma_i$ with the quotient metric, that is,
$$d(a\Gamma_i, b\Gamma_i)=\min\{d(a\gamma_1, b\gamma_2) : \;\; \gamma_1, \gamma_2\in \Gamma_i\}.$$
Let $X(\Gamma)=\bigsqcup_{i=1}^\infty \Gamma/\Gamma_i$ be  the disjoint union of $\Gamma/\Gamma_i$. We
give a metric on $X(\Gamma)$ such that its restriction to each $\Gamma/\Gamma_i$ is the quotient metric defined
above and
$$\lim_{n+m\to \infty, \; n\not= m} d(\Gamma/\Gamma_n, \Gamma/\Gamma_m)=\infty.$$
The metric space $X(\Gamma)$ is called the box metric space \cite{Roe2}.

Recall that the strong Novikov conjecture states that the Baum-Connes  map
$\mu_r: K_{\ast}^{\Gamma}(E\Gamma)\rightarrow K_{\ast}(C^*_r(\Gamma)),$ is injective \cite{K} \cite{BC}, where $E\Gamma$ is the universal
space for free and proper $\Gamma$ actions and $C^*_r(\Gamma)$ is the reduced group $C^*$-algebra.

If $X$ is a discrete metric space with bounded geometry,
the coarse geometric Novikov conjecture states that the Baum-Connes map
$\mu: \lim_{d\rightarrow \infty}K_{\ast}(P_d(X))\rightarrow K_{\ast} (C^*(X)),$ is injective,
where $P_d (X)$ is the Rips complex and $C^*(X)$ is the Roe algebra associated to $X$.
If $X$ doesn't have bounded geometry, then there is a counter-example to the coarse geometric Novikov conjecture \cite{Yu1}.

\begin{thm} If $\Gamma$ has   operator norm localization property and the classifying space $E\Gamma/\Gamma$ for free $\Gamma$-actions has homotopy type of a compact CW
complex,
then the Strong Novikov Conjecture  for $\Gamma$ and all subgroups $\Gamma_n$ ($n=1, 2, 3, \cdots$)
 implies the Coarse Geometric Novikov Conjecture for $X(\Gamma)$.
\end{thm}

Recall that if $\Gamma$ is an infinite property T group, then $X(\Gamma)$ is a sequence of expanders \cite{Lub}. Hence
Theorem 7.1 implies the coarse Novikov conjecture for many interesting examples of sequences of expanders.
A similar result at the level of maximal $C^*$-algebra is proved in \cite{GWY} without the operator norm localization property.

\begin{defn} (Roe  \cite{Roe1}) Let $X$ be a discrete metric space and $T$ be an operator acting on $\ell^2(X)\otimes H$ with finite propagation. $T$ is called locally compact if $T_{x,y}$ is compact for every $x$ and $y$ in $X$, where $T=(T_{x,y})_{x,y\in X}$ is the matrix representation of $T$ with respect to the Hilbert space decomposition $\ell^2(X)\otimes H=\oplus_{x\in X}(\delta_x \otimes H).$
We denote by $\mathbb{C}[X]$ the algebra of all locally compact operators acting on $\ell^2(X)\otimes H$ with finite propagation.
The Roe algebra $C^*(X)$ is the operator norm closure of $\mathbb{C}[X]$.
\end{defn}

If $\Gamma$ is a finitely generated group with a word metric, we denote by $C^*(|\Gamma|)$ the Roe algebra for $\Gamma$ as a metric space with a word metric. If $\Gamma'$ is a subgroup of $\Gamma$, we denote by $\mathbb{C}[|\Gamma| ]^{\Gamma'}$ the fixed
point subalgebra of $\mathbb{C}[|\Gamma| ]$, i.e. $\mathbb{C}[|\Gamma| ]^{\Gamma'}$ consists of  all operators $T$ in $\mathbb{C}[|\Gamma| ]$ satisfying $T_{gx,gy}=T_{x,y}$ for all $g\in \Gamma'$ and $x,y\in \Gamma$. We denote by $C^{\ast}_{r, \Gamma'}(|\Gamma|)$ the operator norm closure of $\mathbb{C}[|\Gamma| ]^{\Gamma'}$.

Let $T\in \mathbb{C}[X(\Gamma)]$. Suppose that $T$ has finite propagation $l$. Let $n$ be the
smallest positive integer such that $d(\gamma, e)>2l$ for all $\gamma\in \Gamma_n$ and
$d_{X(\Gamma)}(\Gamma/\Gamma_i, \Gamma/\Gamma_{j})>2l$ if $i\neq j$ and $i\geq n $ and $j\geq n$, where $e$ is the identity element in
$\Gamma$. Let $$Z=\bigsqcup_{i=1}^{n-1} \Gamma/\Gamma_i,\,\,\,\,\,
Y=\bigsqcup_{i=n}^{\infty} \Gamma/\Gamma_i.$$
$T$ decomposes as follows
$$T=T^0\oplus_{i\geq n} T_i,$$ where
$T^0$ acts on  $\ell^2(Z)\otimes H$ and $T_i$ acts on $\ell^2(\Gamma/\Gamma_i)\otimes H$ for each $i\geq n$. Let
$S_i$ be the operator acting on $\ell^2(\Gamma)\otimes H$  defined by
\[
S_{i;  x, y} = \left\{
\begin{array}{rl}
T_{i; [x], [y]},  & \mbox{if} \; d(x, y)\leq l,     \\
0, &  \mbox{otherwise},
\end{array}
\right.
\]
where, for $x, y\in \Gamma$, $S_{i;  x, y}$ denotes the $(x, y)$-entry of the matrix representation of $S_i$ and,
for $[x], [y]\in \Gamma/\Gamma_i$, the operator $T_{i; [x], [y]}$ is the $([x], [y])$-entry
in the matrix representation  of $T_i$.

We define a map:
$$\phi: \mathbb{C}[X(\Gamma)] \rightarrow \prod_{i=1}^\infty \mathbb{C}[|\Gamma| ]^{\Gamma_i}
\Big/ \bigoplus_{i=1}^\infty \mathbb{C}[|\Gamma| ]^{\Gamma_i}$$
by:
$$\phi(T) =(\oplus_{i<n}0)\oplus \prod_{i\geq n}^\infty S_i. $$
It is not difficult to verify that $\phi$ is a homomorphism.

\begin{lem} If $\Gamma$ has operator norm localization property, then $\phi$ extends to a bounded homomorphism
$$\phi: C^*(X(\Gamma))\rightarrow \prod_{i=1}^\infty C^*_{r, \Gamma_i}(|\Gamma|)
\Big/ \bigoplus_{i=1}^\infty C^*_{r, \Gamma_i}(|\Gamma|).$$
\end{lem}

The proof of this lemma follows from the definition of the operator norm localization property and is therefore
omitted. Now the proof of Theorem 7.1 follows from  our lemma and the argument in the proof of part III of Theorem 5.2 in \cite{GWY}.


\bigskip
\footnotesize

\noindent \noindent Xiaoman Chen\\
Department of Mathematics, Fudan University,\\
 Shanghai 200433, China\\
 E-mail: \url{xchen@fudan.edu.cn}\\

\noindent \noindent Romain Tessera\\
Department of Mathematics, Vanderbilt University,\\ Stevenson
Center, Nashville, TN 37240 United States,\\ E-mail:
\url{romain.a.tessera@vanderbilt.edu}\\

\noindent \noindent Xianjin Wang\\
Department of Mathematics, Fudan University,\\
 Shanghai 200433, China\\
 E-mail: \url{wangxianjin@gmail.com}\\

\noindent \noindent Guoliang Yu\\
Department of Mathematics, Vanderbilt University,\\ Stevenson
Center, Nashville, TN 37240 United States,\\ E-mail:
\url{guoliang.yu@vanderbilt.edu}\\

\end{document}